\documentclass[12pt]{article}

\usepackage[margin=2cm]{geometry}
\usepackage{setspace}
\usepackage{amsmath}
\usepackage{amssymb}
\usepackage{graphicx}
\usepackage{subfig}
\usepackage{algpseudocode}
\usepackage[Algorithm]{algorithm}
\usepackage{tikz}

\usepackage{float}
\usetikzlibrary{patterns}
\usepackage[toc,page]{appendix}

\newcommand{\REQUIRE}{\State{\textbf{Input:\ }}}
\newcommand{\ENSURE}{\State{\textbf{Output:\ }}}
\newcommand{\abs} [1] {\left\lvert{#1}\right\rvert}
\newcommand{\floor} [1] {\left \lfloor{#1}\right \rfloor}
\newcommand{\betless}{\noalign{\vskip3pt plus 3pt minus 1pt}}
\newcommand{\bet}{\noalign{\vskip6pt plus 3pt minus 1pt}}

\renewcommand{\eqref}[1]{\text{Eq.~}\ref{#1}}

\usepackage{color}
\usepackage{xparse}

\newcommand{\rvc}[2]{
\textcolor{red}{#1}%
    \IfNoValueF{#2}{\protect\footnotemark{\footnotetext{\textcolor{red}{#2}}}%
    }%
}

\begin{document}

\title{Direct Inversion of the 3D Pseudo-polar Fourier Transform}
\author
{Amir Averbuch${^1}$~~Gil Shabat${^1}$\footnote{Corresponding author: gil.shabat@cs.tau.ac.il}~~Yoel Shkolnisky${^2}$\\
${^1}$School of Computer Science, Tel Aviv University, Israel\\
${^2}$School of Applied Mathematics, Tel Aviv University, Israel
}
\date{}
\maketitle

\begin{abstract}
	The pseudo-polar Fourier transform is a specialized non-equally spaced Fourier transform, which evaluates the Fourier transform on a near-polar grid known as the pseudo-polar grid. The advantage of the pseudo-polar grid over other non-uniform sampling geometries is that the transformation, which samples  the Fourier transform on the pseudo-polar grid, can be inverted using a fast and stable algorithm. For other sampling geometries, even if the non-equally spaced Fourier transform can  be inverted, the only known algorithms are iterative. The convergence speed of these algorithms as well as their accuracy are difficult to control, as they depend both on the sampling geometry as well as on the unknown reconstructed object.
	In this paper, a direct inversion algorithm for the three-dimensional pseudo-polar Fourier transform is presented. The algorithm is based only on one-dimensional resampling operations, and is shown to be significantly faster than existing iterative inversion algorithms.
\end{abstract}



\section{Introduction}
The fast Fourier transform (FFT) is one of the most commonly used algorithms, and has far reaching implications
in science and technology~\cite{cipra2000best}. While the FFT efficiently computes the discrete Fourier transform (DFT) on a Cartesian grid, in many applications it is required to compute the Fourier transform on non-Cartesian grids. Computing the FFT on a non-Cartesian grid can be implemented using one of the available non-equally spaced FFTs (NUFFTs)~\cite{Dutt_Rokhlin,ware1998,beylkin1995,fessler2003,potts2001,greengard2004}. Of particular importance in applications are the polar grids in two and three dimensions, which emerge naturally in a wide range of applications, from computerized tomography to nano-materials to image processing~\cite{Scott_Chen_Mecklenburg_Zhu_Xu_Ercius_Dahmen_Regan_Miao,Chen_Zhu_White_Chiu_Scott_Regan_Marks_Huang_Miao,Zhao_Brun_Coan_Huang_Sztrokay_Diemoz_Liebhardt_Mittone_Gasilov_Miao_Bravin,Fahimian_Zhao_Huang_Fung_Mao_Zhu_Khatonabadi_DeMarco_Osher_McNitt_Gray_Miao,Jiang_Song_Chen_Xu_Raines_Fahimian_Lu_Lee_Nakashima_Urano_Ishikawa_Tamanoi_Miao,Jiang_Song_Chen_Xu_Murphy_Wright_Jensen_Miao,Lee_Fahimian_Iancu_Suloway_Murphy_Wright_Castano_Diez_Jensen_Miao,kutyniok2012shearlab,bermanis20103,keller2005algebraically}.

Although the Fourier transform can be evaluated on a polar grid by using one of the available NUFFT algorithms, there are two factors that may limit their applicability to large problems. First, although all NUFFT algorithms have the same asymptotic complexity as the classical FFT, namely, $O(n\log n)$, where $n$ is the number of input/output samples, their runtime complexity involves rather larger constants. This is due to local interpolations used by all these methods to resample from the Cartesian grid to the polar grid. Therefore, Although these algorithms are asymptotically as efficient as the FFT, they are prohibitively slow for large input sizes. Second, as many real-life problems can be formulated as the recovery of image samples from frequency samples (e.g., medical
imagery reconstruction algorithms), they actually require computing the inverse Fourier transform. However, unlike the FFT, whose inverse is equal to its adjoint and thus can be easily inverted, the NUFFT, except for very special cases, does not have this property, and therefore, more elaborated inversion schemes are required. In particular, although the transform that evaluates the Fourier transform on the polar grid is formally invertible, the condition number of this transformation excludes such inversion in practice.

Although the polar sampling geometry of Fourier space is natural for various applications such as computerized tomography~\cite{Natterer} or electron microscopy~\cite{Scott_Chen_Mecklenburg_Zhu_Xu_Ercius_Dahmen_Regan_Miao,Chen_Zhu_White_Chiu_Scott_Regan_Marks_Huang_Miao}, it is sometimes possible to overcome the aforementioned difficulties by replacing the polar frequency sampling grid with the pseudo-polar frequency sampling grid. The pseudo-polar grid~\cite{averbuch2008framework} (to be described below) is nearly polar, but unlike the polar grid, it is not equally spaced in the angular direction, but uses equally spaced slopes. We refer to the transformation which samples the Fourier transform on the pseudo-polar grid, as the pseudo-polar Fourier transform (PPFT). Both the PPFT and its inverse admit efficient and stable algorithms.  Iterative algorithms used to invert the pseudo-polar Fourier transforms can be found in~\cite{averbuch2008framework,grochenig1992reconstruction,feichtinger1995efficient,rauth1998smooth}. Since the pseudo-polar Fourier transform is ill-conditioned, these inversion algorithms require preconditioning. Several preconditioning approaches have been suggested, such as using Voronoi weights~\cite{fenn2007computation} or a preconditioner derived from the Jacobian~\cite{iterativeradon2014}. In particular, it is shown in~\cite{iterativeradon2014} that under the appropriate preconditioning, the condition number of the pseudo-polar Fourier transform is small. Therefore, it can be efficiently inverted using iterative inversion schemes. Nevertheless, such an iterative inversion requires to compute the pseudo-polar Fourier transform and its adjoint in each iteration, which as mentioned above, is too expensive for large problems. Even if the pseudo-polar Fourier transform and its adjoint are computed using specialized algorithms~\cite{averbuch20033d}, the number of iterations required for the inversion may be of the order of a few dozens, which may be too slow for large scale problems. Moreover, the accuracy of the reconstructed object by such iterative methods typically cannot be estimated from the convergence criterion and more calculations are needed to verify that the reconstructed object has the required accuracy.

The PPFT can also be expressed as a multilevel Toeplitz operator \cite{rauth1998smooth}. Thus, the inverse PPFT is equivalent to the inversion of this operator. Direct and iterative algorithms have been devised for inverting Toeplitz-related operators. Direct inversion algorithms, which are stable and fast, are available for Toeplitz matrices inversion~\cite{stewart2003superfast,chandrasekaran2007superfast}. On the other hand, these algorithms do not extend to multilevel-Toeplitz operators. The only available exact algorithm for inverting two-level Toeplitz operators~\cite{turnes2012image} is not applicable to our setting since it assumes that the two-level Toeplitz matrices are triangular in one or both their levels. 

In this paper, we present a direct algorithm for inverting the three-dimensional pseudo-polar Fourier transform, which generalizes the two-dimensional direct inversion algorithm in~\cite{averbuch2008framework}. Unlike the iterative approaches, the proposed algorithm terminates within a fixed number of operations, which is determined only by the size of the input. The inversion algorithm is shown numerically to have high accuracy, and since most of its steps are implemented in-place, its memory requirements are sufficiently low to be applicable to very large inputs. Finally, the algorithm is numerically stable since it consists of a series of well-conditioned steps. In particular, no preconditioner is needed.

This paper is organized as follows. In Section~\ref{pseudo-polar-background}, we provide the required background on the pseudo-polar grid and the pseudo-polar Fourier transform. In Section~\ref{sec:onion}, we give a detailed description of
the direct inversion algorithm for the three-dimensional pseudo-polar transform. Numerical results and comparisons to other algorithms are given in
Section~\ref{numerical}. Finally, some concluding remarks are given in Section~\ref{sec:conclusions}.

\section{Mathematical preliminaries}
\label{pseudo-polar-background}

In this section, we provide the mathematical background required to derive our inversion algorithm. In Section~\ref{sec:pseudopolarft} we revisit the pseudo-polar Fourier transform. In Section~\ref{sec:2D Radon:direct inverse solving toeplitz systems}, we describe a fast algorithm for solving Toepliz systems of equations. This algorithm is used in Section~\ref{sec:1d resampling algorithm} for fast resampling of univariate trigonometric polynomials.

\subsection{Pseudo-polar Fourier transform}
\label{sec:pseudopolarft}
\begin{figure}
	\centering
	\subfloat[The sector $\Omega_{ppx}$ for $n=4$]{\label{ppx}
		\includegraphics[height=4.5cm]{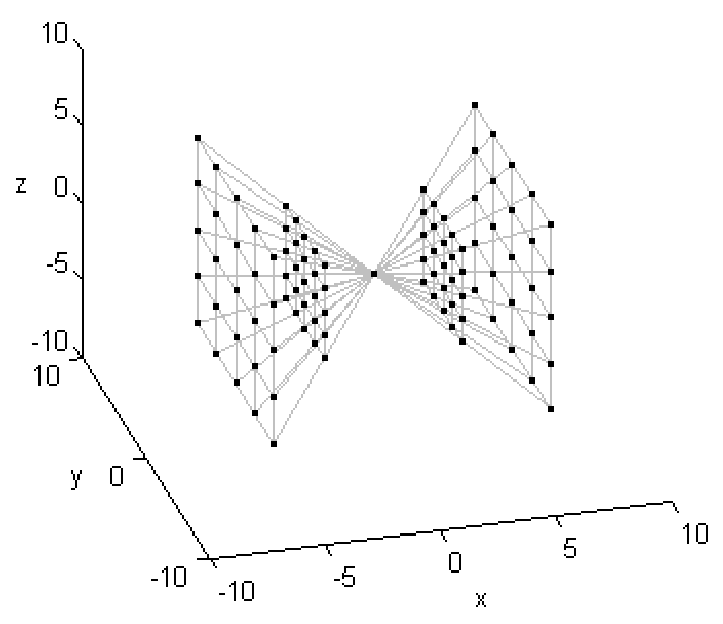}}
	\subfloat[The sector $\Omega_{ppy}$ for $n=4$]{\label{ppy}
		\includegraphics[height=4.5cm]{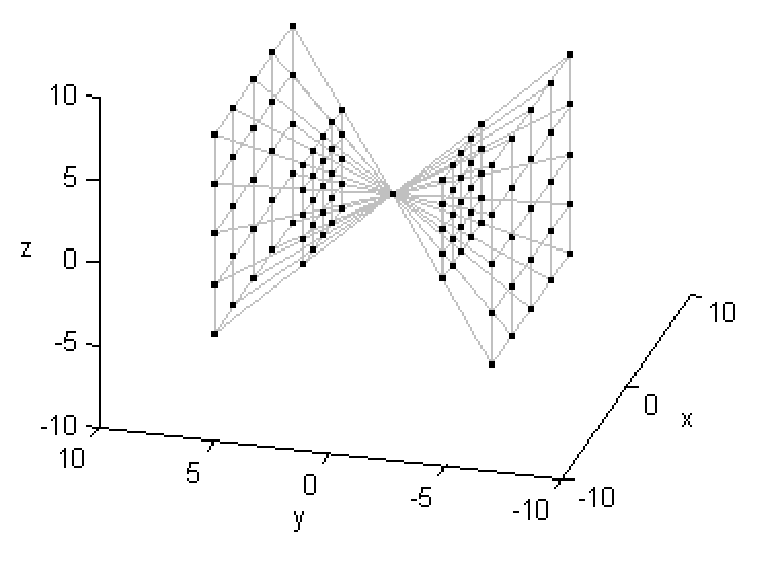}}
	\subfloat[The sector $\Omega_{ppz}$ for $n=4$]{\label{ppz}
		\includegraphics[height=4.5cm]{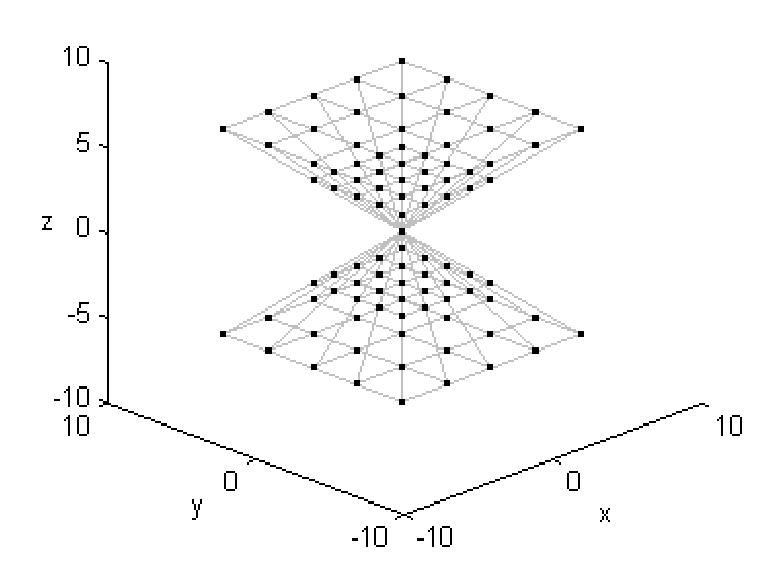}}
	\caption{The three sectors $\Omega_{ppx},\Omega_{ppx}$ and $\Omega_{ppz}$ that form the 3D pseudo-polar grid.}
	\label{fig:ppgrid}
\end{figure}

The 3D pseudo-polar grid, denoted $\Omega_{pp}$, is
defined by
\begin{equation}\label{eq:ppgrid}
\Omega_{pp}=\Omega_{ppx} \cup \Omega_{ppy} \cup
\Omega_{ppz},
\end{equation}
where
\begin{equation}
\label{eq2.4}
\begin{array}{lll}
\Omega_{ppx}&= & \left\{\left(k,-\frac{2l}{n}k,-\frac{2j}{n}k\right):l,j=-n/2,\ldots,n/2 , k=-\floor{\frac{m}{2}},\ldots,\floor{\frac{m}{2}} \right\},
\\
\Omega_{ppy}&= & \left\{\left(-\frac{2l}{n}k,k,-\frac{2j}{n}k\right):l,j=-n/2,\ldots,n/2 , k=-\floor{\frac{m}{2}},\ldots,\floor{\frac{m}{2}} \right\},
\\
\Omega_{ppz}&= &
\left\{\left(-\frac{2l}{n}k,-\frac{2j}{n}k,k\right):l,j=-n/2,\ldots,n/2 , k=-\floor{\frac{m}{2}},\ldots,\floor{\frac{m}{2}} \right\},
\end{array}
\end{equation}
with $m=qn+1$ for an even $n$ and a positive integer $q$. We denote a specific point in $\Omega_{ppx}$, $\Omega_{ppy}$ and $\Omega_{ppz}$
by $\Omega_{ppx}(k,l,j)$, $\Omega_{ppy}(k,l,j)$ and
$\Omega_{ppz}(k,l,j)$, respectively. The psuedo-polar grid is illustrated in Fig.~\ref{fig:ppgrid} for $n=4$ and $q=1$ . As can be seen from Fig.~\ref{fig:ppgrid}, the pseudo-polar grid consists of equally spaced samples along rays,
where different rays have equally spaced slopes but the angles between adjacent rays are not equal. This is the key difference between the pseudo-polar grid and the polar grid~\cite{averbuch2006fast}. Thus, we can refer to $k$ as a ``pseudo-radius" and to $l$ and $j$ as ``pseudo-angles".

Next, we define the discrete time Fourier transform of an $n \times n \times n$ volume $I$ by
\begin{equation}
\label{eq:3DPPFT:dtft3}
\hat{I}(\omega_{x},\omega_{y},\omega_{z})=\sum_{u,v,w=-n/2}^{n/2-1} I(u,v,w) \, e^{2\pi \imath (u\omega_{x}+v\omega_{y}+w\omega_{z})/m}, \quad \omega_{x},\omega_{y},\omega_{z} \in [-m/2,m/2],
\end{equation}
where as before $m=qn+1$ for an even $n$ and a positive integer $q$. The three-dimensional pseudo-polar Fourier transform (PPFT) is defined as the samples of the discrete time Fourier transform of~\eqref{eq:3DPPFT:dtft3} on the pseudo-polar grid $\Omega_{pp}$. Specifically, if we denote the concatenation of three arrays $A_{1}$, $A_{2}$, and $A_{3}$ by
$A=\left[A_{1}, A_{2}, A_{3} \right ]$, then the pseudo-polar Fourier transform of a volume $I \in \mathbb{C}^{n\times n\times  n}$ is an array $\hat{I}_{\Omega_{pp}}\in\mathbb{C}^{3\times(n+1)\times(n+1)\times(3n+1)}$ given by 
\begin{equation}\label{eq:ppft}
\hat{I}_{\Omega_{pp}}= \left [ \hat{I}_{\Omega_{ppx}},\hat{I}_{\Omega_{ppy}},\hat{I}_{\Omega_{ppz}}\right ],
\end{equation} 
where
\begin{equation}\label{eq:Omega2}
\hat{I}_{\Omega_{ppx}}=\hat{I}\left(k,-\frac{2l}{n}k,-\frac{2j}{n}k\right), \quad \hat{I}_{\Omega_{ppy}}=\hat{I}\left(-\frac{2l}{n}k,k,-\frac{2j}{n}k\right), \quad
\hat{I}_{\Omega_{ppz}}=\hat{I}\left(-\frac{2l}{n}k,-\frac{2j}{n}k,k\right),
\end{equation}
and $\Omega_{ppx}, \Omega_{ppy}$ and $\Omega_{ppz}$ are defined in~\eqref{eq2.4}. The parameter $q$ in Eqs.~\ref{eq:ppgrid} and ~\ref{eq:3DPPFT:dtft3} determines the frequency resolution of the transform, as well as its geometric properties. For example, as shown in~\cite{averbuch20033d}, in order to derive a three-dimensional discrete Radon transform based on the pseudo-polar Fourier
transform, $q$ must satisfies $q \ge 3$.

The pseudo-polar grid appeared in the literature
several times under different names. It was originally introduced by~\cite{pasciak1973note} under the name
``Concentric Squared Grid" in the context of computerized tomography. More recent works in the context of computerized tomography, which take advantage of the favorable numerical and computational properties of the grid include~\cite{mao2010development,miao2005equally}, where equally sloped tomography is used for radiation dose reduction. Other image processing applications that use the pseudo-polar grid include synthetic aperture radar (SAR) imaging~\cite{lawton1988new}, Shearlets~\cite{kutyniok2012shearlab,kutyniok2012digital},  registration~\cite{liu2006pseudo}, and denoising~\cite{starck2002curvelet}, to name a few.

Recently,~\cite{iterativeradon2014} proposed fast iterative inversion algorithms for the 2D and 3D
pseudo-polar Fourier transforms, which are based on the well-known convolution structure of their Gram operators, combined with a preconditioned conjugate gradients~\cite{golub2012matrix} solver. The convolution structure of the transform allows to invert it using the highly
optimized FFT algorithm instead of the forward and adjoint
transforms derived in~\cite{averbuch2008framework}. As the conjugate gradients iterations require preconditioning,~\cite{iterativeradon2014} proposes a preconditioner that leads to a
very small condition number, thus making the inversion process fast and accurate.
However, this approach has drawbacks such as high memory
requirements and dependency of the number of iterations and on the size of the data.
The algorithm presented in this paper, which has a low memory requirements, achieves high accuracy and is faster than the iterative algorithm~\cite{iterativeradon2014}.

\subsection{Solving Toeplitz systems}
\label{sec:2D Radon:direct inverse solving toeplitz systems}

Let $A_n$ be an $n \times n$ Toeplitz matrix and let $y$
be an arbitrary vector of length $n$. We describe a fast
algorithm for computing $A_n^{-1}y$. This algorithm is well-known and appears, for example, in~\cite{averbuch2008framework}, but we repeat it here for completeness of the description. The algorithm  consists of a fast factorization of the inverse Toeplitz matrix followed by a fast
algorithm that applies the inverse matrix to a vector~\cite{Gohberg1994,Kailath99}. We denote by $T_{n}(c,r)$ an $n \times n$ Toeplitz
matrix whose first column and row are $c$ and $r$, respectively. For symmetric matrices, $c = r$.

Circulant matrices are diagonalized by the Fourier matrix. Hence,
a circulant matrix $C_n$ can be written as $C_{n} = W^{*}_{n} D_{n}
W_{n}, $ where $D_{n}$ is a diagonal matrix containing the
eigenvalues $\lambda_{1},\ldots,\lambda_{n}$ of $C_{n}$, and
$W_{n}$ is the Fourier matrix given by $W_{n}(j,k)=
\frac{1}{\sqrt{n}} e^{2 \pi \imath j k /n}$. Moreover, if
$c=[c_{0},c_{1},\ldots,c_{n-1}]^{T}$ is the first column of
$C_{n}$, then $W_{n}c=[\lambda_{1},\ldots,\lambda_{n}]^{T}$.
Obviously, the matrices $W_{n}$ and $W_{n}^{*}$ can be applied in
$O(n \log{n})$ operations using the FFT. Thus, the
multiplication of $C_{n}$ with an arbitrary vector $x$ of length
$n$ can be implemented in $O(n \log{n})$ operations by applying
FFT to $x$, multiplying the result by $D_{n}$, and then applying the
inverse FFT.

To compute $A_{n}x$ for an arbitrary Toeplitz matrix
$A_{n}=T_{n}(c,r)$ and an arbitrary vector $x$, we first embed
$A_{n}$ in a $2n \times 2n$ circulant matrix $C_{2n}$
\begin{equation*}
C_{2n} = \left (
\begin{array}{cc}
A_{n} & B_{n} \\
\betless
B_{n} & A_{n}
\end{array}
\right ),
\end{equation*}
where $B_{n}$ is an $n \times n$ Toeplitz matrix given by
\begin{equation*}
B_{n}=T_{n}([0,r_{n-1},\ldots,r_{2},r_{1}],[0,c_{n-1},\ldots,c_{2},c_{1}]).
\end{equation*}
Then, $A_{n} x$ is computed in $O(n \log{n})$ operations by zero
padding $x$ to length $2n$, applying $C_{2n}$ to the padded
vector, and discarding the last $n$ elements of the resulting vector.

Next, assume that $A_{n}$ is invertible. The Gohberg--Semencul
formula~\cite{Gohberg1994,Gohberg1972} provides an explicit representation
of $A_{n}^{-1}$ as
\begin{equation}\label{eq:gohberg-semencul formula}
A_{n}^{-1} = \frac{1}{x_{0}} \left ( M_{1} M_{2} -M_{3} M_{4}
\right ),
\end{equation}
where
\begin{align}
M_{1} &= T_{n}([x_{0},x_{1},\ldots,x_{n-1}],[x_{0},0,\ldots,0]),  \label{ppft:m1eq} \\
\bet
M_{2} &= T_{n}([y_{n-1},0,\ldots,0],[y_{n-1},y_{n-2},\ldots,y_{0}]), \label{ppft:m2eq} \\
\bet
M_{3} &= T_{n}([0,y_{0},\ldots,y_{n-2}],[0,\ldots,0]), \label{ppft:m3eq} \\
\bet
M_{4} &= T_{n}([0,\ldots,0],[0,x_{n-1},\ldots,x_{1}]), \label{ppft:m4eq}.
\end{align}
The vectors $x=[x_{0},\ldots,x_{n-1}]$ and $y=[y_{0},\ldots,y_{n-1}]$ are given as the solutions to
\begin{equation}\label{eq:xyequation}
\begin{alignedat}{5}
A_{n} x &=& e_{0},   &\qquad &e_{0}&=[1,0,\ldots,0]^{T}, \\
A_{n} y &=& e_{n-1}, &\qquad &e_{n-1}&=[0,\ldots,0,1]^{T}.
\end{alignedat}
\end{equation}

The matrices $M_{1}$, $M_{2}$, $M_{3}$, and $M_{4}$ have Toeplitz
structure and are represented implicitly using the vectors $x$ and
$y$. Hence, the total storage required to store $M_{1}$, $M_{2}$,
$M_{3}$, and $M_{4}$ is that of $2n$ numbers. If the matrix $A_{n}$ is
fixed, then the vectors $x$ and $y$ can be precomputed. Once the
triangular Toeplitz matrices $M_{1}$, $M_{2}$, $M_{3}$ and $M_{4}$ have been computed, the application of $A_{n}^{-1}$ is reduced to the
application of four Toeplitz matrices. Thus, the application of
$A_{n}^{-1}$ to a vector requires $O(n \log{n})$ operations. The pseudo-code of applying $A_{n}^{-1}$ to a vector is described in Algorithms~\ref{alg:3DPPFT:topdiag},~\ref{alg:3DPPFT:topinvmul}, and~\ref{alg:3DPPFT:topmul} in the appendix.
Algorithm~\ref{alg:3DPPFT:topdiag} lists the function \texttt{ToeplitzDiag}, which computes the diagonal form of the circulant embedding of a Toeplitz matrix. Algorithm~\ref{alg:3DPPFT:topinvmul} lists the function \texttt{ToeplitzInvMul}, which efficiently multiplies an inverse Toeplitz matrix, given in its diagonal form, by a vector. The latter function uses \texttt{ToeplitzMul}, listed in Algorithm~\ref{alg:3DPPFT:topmul}, which efficiently multiplies a general Toeplitz matrix by a vector.

\subsection{Resampling trigonometric polynomials}\label{sec:1d resampling algorithm}
The main (and in fact the only) tool behind our algorithm is resampling of univariate trigonometric polynomials. Assume we are given a set of points $\{ y_j \}_{j=1}^{N}$, $y_j \in [-\pi,\pi]$, and their values $\{ f(y_j) \}_{j=1}^{N}$, where $f$ is some unknown univariate trigonometric polynomial of degree $n\leq N$. We want to estimate the values of $f$ at a new set of points $\{ x_j \}_{j=1}^{M}$, $x_j \in [-\pi,\pi]$. This can be formulated by first solving
\begin{equation}\label{eq:resamplingLS}
\min_{\alpha} \sum_{j=1}^{N} \left \vert f(y_j)-\sum_{k=-n/2}^{n/2-1} \alpha_k e^{\imath k y_j} \right \vert^2, \quad \alpha=\left (a_{-n/2},\ldots,\alpha_{n/2-1} \right ) \in \mathbb{C}^{n},
\end{equation}
for the coefficients vector $\alpha$, followed by evaluating
\begin{equation}\label{eq:fxi}
f(x_{i}) = \sum_{k=-n/2}^{n/2-1} \alpha_k e^{\imath k x_i }.
\end{equation}
In matrix notation,~\eqref{eq:resamplingLS} is written
\begin{equation}
\label{eq:ls_ppft}
\min_{\alpha} \Vert f-A\alpha \Vert_2,
\end{equation}
where the entries of the matrix $A$ are given by $a_{jk}=e^{\imath k y_j }$ and the coordinates of the vector $f$ are $f_j=f(y_j)$ (we denote by $f$ both the vector and the function, as the appropriate meaning is clear from the context).
Direct solution for the coefficients vector $\alpha$ requires $\mathcal{O}(n^3)$ operations (assuming $N$ and $M$ are of size $O(n)$). Obviously, solving for $\alpha$ also depends on the condition number of $A$, which affects the accuracy of $f(x_i)$.

We next present a fast algorithm for computing $f(x_{i})$ (faster than directly solving first for $\alpha$), which exploits the Toeplitz structure of $A^{*}A$ and uses the non-equally spaced FFT (NUFFT)~\cite{Dutt_Rokhlin,ware1998,beylkin1995,fessler2003,potts2001,greengard2004}. The resulting algorithm has complexity of  $\mathcal{O}(n \log n + n \log 1/\epsilon)$ operations, where $\epsilon$ is the accuracy of the computations, in addition to a preprocessing step that takes $\mathcal{O}(n^2)$ operations.

Following~\cite{greengard2004accelerating}, we define the NUFFT type-I by
\begin{equation}
f_k=\frac{1}{n_j}\sum_{j=1}^{n_j}c_j e^{\pm \imath k x_j}, \quad k=-n/2,\ldots,n/2-1.
\label{eq:nuff1}
\end{equation}
and the NUFFT type-II by
\begin{equation}
c_j=\sum_{k=-n/2}^{n/2-1} f_k e^{\pm \imath k x_j}, \quad j=1,...,n_j.
\label{eq:nuff2}
\end{equation}
Both types can be approximated to a relative accuracy $\epsilon$ in $\mathcal{O}(n \text{log}n + n \log 1/\epsilon)$ operations by any of the aforementioned NUFFT algorithms. From the definitions of NUFFT of types I and II, we see that $A^*f$, where $A$ is the matrix from~\eqref{eq:ls_ppft} and $f$ is an arbitrary vector, is equal to the application of NUFFT type-I to $f$. Similarly, the application of $A$ to an arbitrary vector $c$ is equal to the application of NUFFT type-II to $c$. Hence, the application of $A$ or $A^*$ to a vector can be implemented in $\mathcal{O}(n \log n + n \log 1/\epsilon)$ operations.

To solve the least-squares problem of~\eqref{eq:ls_ppft} we form the normal equations
\begin{equation}
\label{eq:ppft:LSAstar}
A^*A\alpha =A^*f.
\end{equation}
The right-hand side in~\eqref{eq:ppft:LSAstar} can be computed efficiently using NUFFT type-I. The matrix $A^*A$ is a symmetric Toeplitz matrix of size $n\times n$. The first column of $A^*A$, which due to the symmetric Toeplitz structure encodes all its entries, is computed efficiently by applying NUFFT type-I to the vector $w$ whose entries are $w_j=e^{-\imath n y_j/2}$.
Computing $(A^*A)^{-1}$  takes $\mathcal{O}(n^2)$ operations using the Durbin-Levinson algorithm~\cite{levinson1947wiener} and applying it to a vector takes $\mathcal{O}(n\log n)$ operations using the Gohberg-Semencul formula~\cite{kailath1989generalized}, as was described in Section~\ref{sec:2D Radon:direct inverse solving toeplitz systems}. This procedure is described in details in~\cite{averbuch2008framework}.
Since $(A^*A)^{-1}$ depends only on $n$, it can be  precomputed. Therefore, solving for $\alpha$ in~\eqref{eq:ppft:LSAstar} takes $\mathcal{O} (n \log n + n \log 1/\epsilon)$ operations and computing $f(x_{i})$ in~\eqref{eq:fxi} takes additional  $\mathcal{O} (n\log n+n\log 1/\epsilon)$ operations using NUFFT type-II. The entire resampling algorithm is described in Algorithm~\ref{alg:3DPPFT:toepresampling}.

\section{Direct inversion of the 3D PPFT}
\label{sec:onion}

Given the PPFT $\hat{I}_{\Omega_{pp}}$ defined in~\eqref{eq:ppft} of an unknown $n\times n \times n$ volume $I$, the proposed direct inversion algorithm recovers $I$ in two steps. The first
step resamples the PPFT into an intermediate Cartesian
grid. Specifically, we resample the trigonometric polynomial $\hat{I}$ of~\eqref{eq:3DPPFT:dtft3} from the grid $\Omega_{pp}$ in~\eqref{eq:ppgrid} to the frequency grid
\begin{equation}
\label{eq:DecCaGrid}
\Omega_{c}=\{(qu,qv,qw): u,v,w=-n/2,...,n/2 \} .
\end{equation}
The second step of our algorithm recovers $I$ from the samples of $\hat{I}$ on $\Omega_{c}$. Note that this second step cannot be directly implemented by inverse FFT. The two steps of the algorithm are described in Sections~\ref{sec:onion-peeling} and~\ref{sec:invDecimatedF}, respectively, its pseudo-code is given in Algorithm~\ref{alg:3DPPFT:inversion}, and its complexity is analyzed in Section~\ref{sec:complexity}.

\subsection{Resampling the pseudo-polar grid to the grid $\Omega_{c}$}\label{sec:onion-peeling}

We start by describing the procedure for resampling $\hat{I}$ from $\Omega_{pp}$ to $\Omega_{c}$. It is based on an ``onion peeling" approach, which resamples at iteration $k$, $k=n/2,\ldots,0$,  points in $\Omega_{pp}$ with ``pseudo-radius'' $k$ to points in $\Omega_{c}$ that lie on a plane. When the iteration that corresponds to $k=0$ is completed,  the values of $\hat{I}$ have been computed on all points of the grid $\Omega_{c}$. We define the sets $\Omega_{c}^{(k)} \subseteq \Omega_{c}$, $k=n/2,\ldots,0$, consisting of the frequencies on which $\hat{I}$ has been resampled until (and including) iteration $k$. Note that for convenience, we index the iterations from $k=n/2$ to $k=0$. Thus, we have that $\Omega_{c}^{(k)} \subseteq \Omega_{c}^{(k-1)}$ and $\Omega_{c}^{(0)}=\Omega_{c}$. Formally, we define
$\Theta^{(k)}$ to be the frequencies on which we resample $\hat{I}$ at iteration $k$ so that
\begin{equation*}
\Omega_{c}^{(k)} = \Omega_{c}^{(k+1)} \cup \Theta^{(k)},
\end{equation*}
where $\Theta^{(k)}$ is give by
\begin{equation}
\label{eq:thetak}
\Theta^{(k)} = \Theta^{(k)}_{+x} \cup \Theta^{(k)}_{-x} \cup \Theta^{(k)}_{+y} \cup \Theta^{(k)}_{-y}
\cup \Theta^{(k)}_{+z}\cup \Theta^{(k)}_{-z} ,
\end{equation}
and for $j,l=-k,\ldots,k$
\begin{equation} 
\begin{alignedat}{3}
\Theta^{(k)}_{+x} &=&  \{(qk,qj,ql)\}, \quad &\Theta^{(k)}_{-x} &=&  \{(-qk,qj,ql)\},\\
\Theta^{(k)}_{+y} &=&  \{(qj,qk,ql)\}, \quad &\Theta^{(k)}_{-y} &=&  \{(qj,-qk,ql)\},\\
\Theta^{(k)}_{+z} &=&  \{(qj,ql,qk)\}, \quad &\Theta^{(k)}_{-z} &=&  \{(qj,ql,-qk)\}.
\end{alignedat}
\label{eq:thetagrid}
\end{equation}
For $k=n/2$ we set $\Omega_{c}^{(n/2)} = \Theta^{(n/2)}$. Moreover, the frequencies of $\Theta^{(n/2)}$ are the frequencies $\Omega_{ppx}(\pm n/2,j,l)$, $\Omega_{ppy}(\pm n/2,j,l)$, and $\Omega_{ppz}(\pm n/2,j,l)$, $l,j=-n/2,\ldots,n/2$, defined in Eq.~\ref{eq2.4}. Thus, the values of $\hat{I}$ on $\Omega_{c}^{(n/2)}$ are given as a subset of the values of $\hat{I}$ on $\Omega_{pp}$. Therefore, no resampling is needed at this iteration.

For subsequent iterations, we use an example to accompany and clarify the formal description. We use as an example a grid corresponding to $n=8$ that demonstrates how the values of $\hat{I}$ on $\Theta^{(n/2-1)}$ (which is $\Theta^{(3)}$ in our particular example) are recovered. Since the same procedure recovers the values of $\hat{I}$ on $\Theta^{(k)}_{\pm x}$, $\Theta^{(k)}_{\pm y}$, and $\Theta^{(k)}_{\pm z}$, we only explain the process for the grid $\Theta^{(k)}_{+x}$.  The pseudo-code for recovering the values of $\hat{I}$ on $\Theta^{(k)}$ is given in Algorithm~\ref{alg:3DPPFT:recover2d}, whose details are also explained below.

\begin{figure}
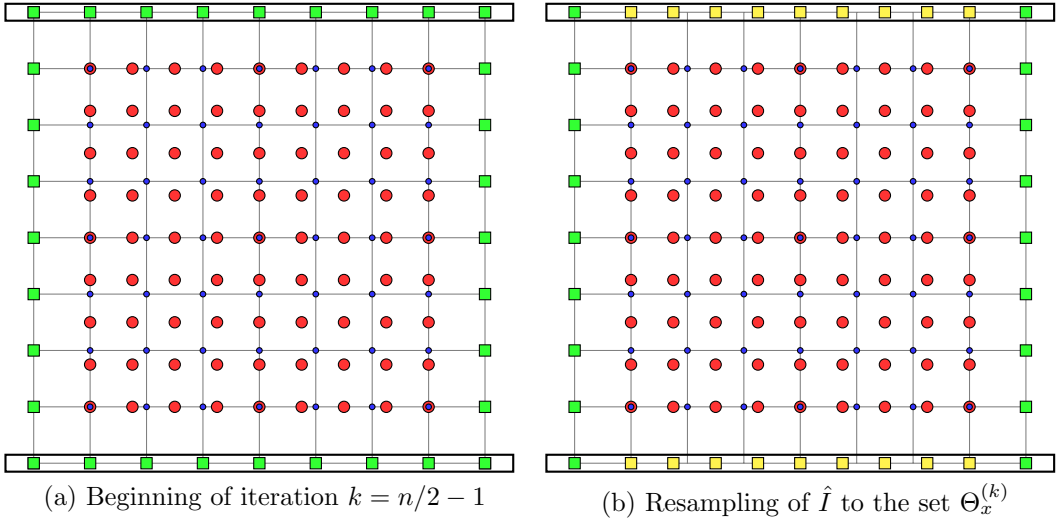

	\begin{center}
		\subfloat[Beginning of iteration $k=n/2-1$]{
			\input{figstep1a.tkz}
			\label{fig:resampling_a_1}
		}
		\subfloat[Resampling of $\hat{I}$ to the set $\Theta^{(k)}_{x}$]{
			\input{figstep1b.tkz}
			\label{fig:resampling_a_2}
		}
	\end{center}
	\caption{The grid before and after the applications of lines $11-14$ in Algorithm~\ref{alg:3DPPFT:recover2d}. Green solid squares correspond to points from $\Omega_{c}^{(k+1)}$, red circles correspond to points from $\Omega_{ppx}$ with a fixed $k=n/2-1$, blue dots correspond to points of $\Theta^{(k)}_{+x}$, and patterned yellow points correspond to point of $\Omega_{c}^{(k+1)}$ resampled to their intermediate position.} \label{fig:resampling_a}
\end{figure}

In Fig.~\ref{fig:resampling_a_1}, green solid squares correspond to points from $\Omega_{c}^{(k+1)}$ on which $\hat{I}$ has already been evaluated (in previous iterations), red circles correspond to points from $\Omega_{ppx}$ with a fixed $k=n/2-1$, and blue dots correspond to points of $\Theta^{(k)}_{+x}$ on which we want to evaluate $\hat{I}$. The frequencies in $\Theta^{(k)}_{+x}$ are points in $\mathbb{R}^{3}$, however, as can be seen from~\eqref{eq:thetagrid}, they all lie on the same plane, and therefore, we depict them as a two-dimensional image whose axes are $i$ and $j$ from~\eqref{eq:thetagrid}.

\begin{figure}
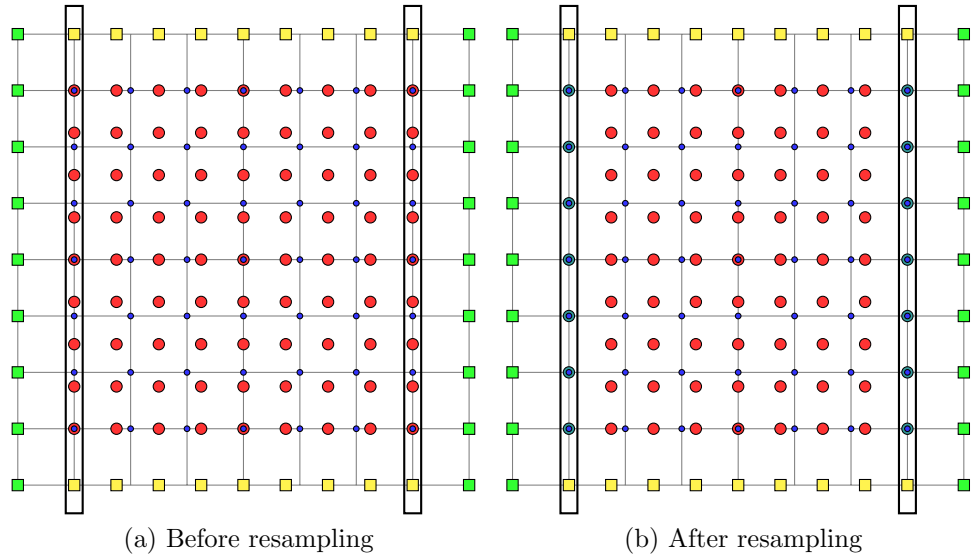

	\begin{center}
		\subfloat[Before resampling]{
			\input{figstep2a.tkz}
			\label{fig:resampling_b_1}
		}
		\subfloat[After resampling]{
			\input{figstep2b.tkz}
			\label{fig:resampling_b_2}
		}
	\end{center}
	\caption{The grid before and after the applications one iteration of lines $18-20$ in Algorithm~\ref{alg:3DPPFT:recover2d}.}
	\label{fig:resampling_b}
\end{figure}

\begin{figure}
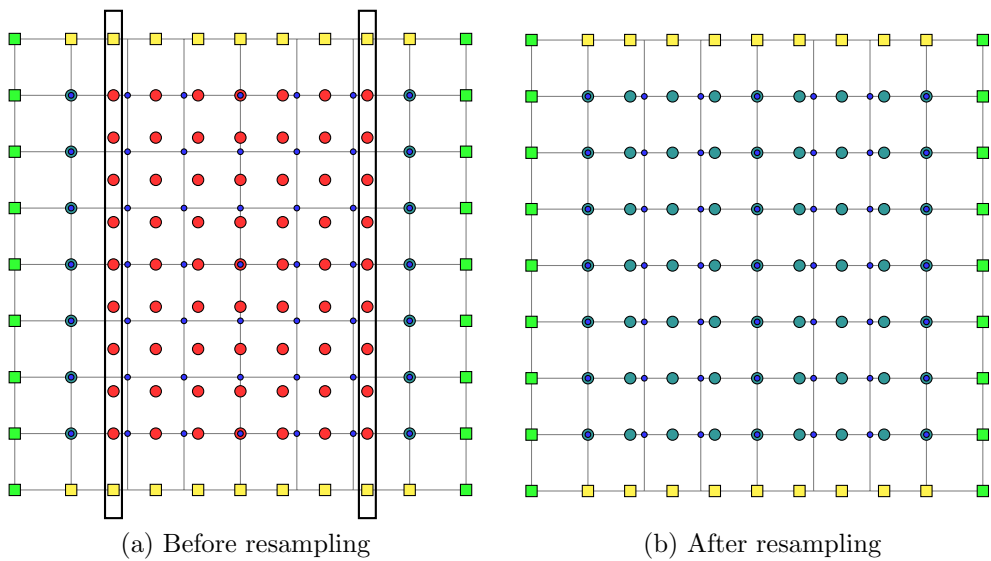

	\begin{center}
		\subfloat[Before resampling]{
			\input{figstep2c.tkz}
			\label{fig:resampling_b_3}
		}
		\subfloat[After resampling]{
			\input{figstep2d.tkz}
			\label{fig:resampling_b_4}
		}
	\end{center}
	\caption{The grid before and after further iterations of lines $18-20$ in Algorithm~\ref{alg:3DPPFT:recover2d}.} \label{fig:gridj1_resampled2}
\end{figure}

\begin{figure}
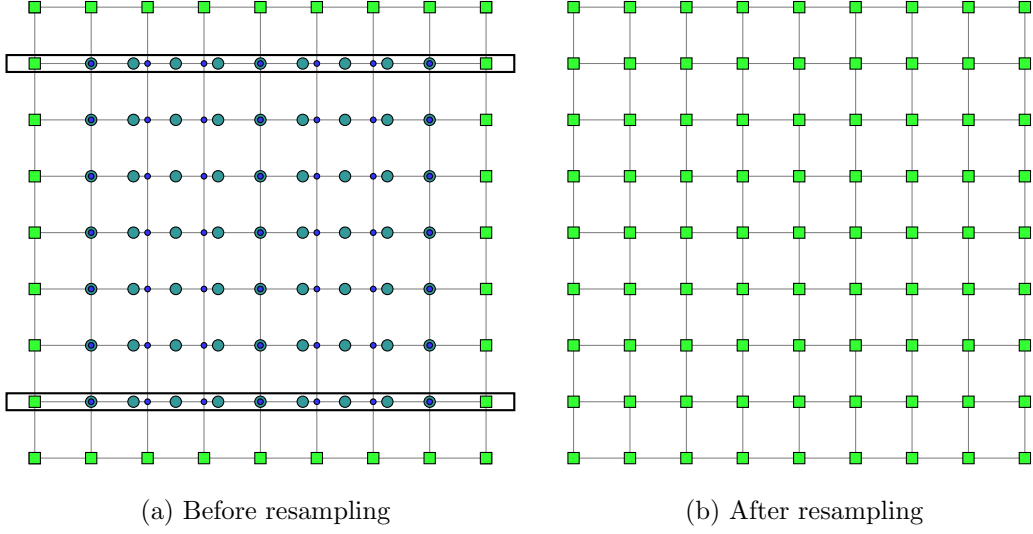

	\begin{center}
		\subfloat[Before resampling]{
			\input{figstep3a.tkz}
			\label{fig:resampling_c_1}
		}
		\subfloat[After resampling]{
			\input{figstep3b.tkz}
			\label{fig:resampling_c_2}
		}
		\caption{The grid before and after the application of lines $22-24$ in Algorithm~\ref{alg:3DPPFT:recover2d} (row resampling).}
	\end{center}
\end{figure}

The first step, depicted in Fig.~\ref{fig:resampling_a_2}, consists of resampling the points of $\Omega_{c}^{(k+1)}$ (solid green squares) to the same spacing as the pseudo-polar points. In the case of $k=n/2-1$, this means that the first and the last rows are resampled, and the result of this resampling is denoted by patterned yellow squares. This step is implemented by lines $11-14$ in Algorithm~\ref{alg:3DPPFT:recover2d}.  Next, as depicted in Fig.~\ref{fig:resampling_b_1}, for all columns, we use the latter resampled points (patterned yellow squares) together with the points of the pseudo-polar grid at the same column to resample $\hat{I}$ to intermediate sampling points. Figures~\ref{fig:resampling_b_1} and~\ref{fig:resampling_b_2} depict the resampling of one such column. Note that the resampled values of $\hat{I}$ after this step, which are depicted as filled teal circles in Fig.~\ref{fig:resampling_b_2}, are not yet the values of $\hat{I}$ on the points of $\Theta^{(k)}_{x}$. Up to here, we only applied resampling to the columns (but not yet to rows). We repeat this for all the columns (see Fig.~\ref{fig:resampling_b_3} for resampling of another column), which results in the grid of Fig.~\ref{fig:resampling_b_4}. This is implemented by lines $18-20$ of Algorithm~\ref{alg:3DPPFT:recover2d}, which using the conventions of Fig.~\ref{fig:gridj1_resampled2}, resamples the red circles along the columns to the filled teal circle, which are now on the same rows as the points indicated by the blue dots (our target sampling points). Next, we apply one-dimensional resampling to each row, by using the resampled points from the previous steps together with the points of $\Omega_{c}^{(k+1)}$ (Fig.~\ref{fig:resampling_c_1}) to recover the samples of $\Theta_{x}^{(k)}$  (Fig.~\ref{fig:resampling_c_2}). This is implemented by lines $22-24$ of Algorithm~\ref{alg:3DPPFT:recover2d}. All one-dimensional resampling operations in Algorithm~\ref{alg:3DPPFT:recover2d} are implemented using the function \texttt{ToeplitzResample} in Algorithm~\ref{alg:3DPPFT:toepresampling}.

\begin{algorithm}
	\caption{Direct inversion of the 3D pseudo-polar Fourier transform.}
	\begin{algorithmic}[1]
		\REQUIRE Samples of the 3D pseudo-polar Fourier transform
		$\hat{I}_{\Omega_{ppx}}, \hat{I}_{\Omega_{ppy}},
		\hat{I}_{\Omega_{ppz}}$, given by~\eqref{eq:Omega2}, each of size  $(qn+1)\times(n+1)\times(n+1)$.
		\ENSURE Volume $I$ of size $n \times n \times n$.		
		\State{ $\hat{I}_D \leftarrow \operatorname{zeros}(n+1,n+1,n+1)$}
		\For{$j=1,\ldots,n/2$}
		
		\State{$\hat{I}_D(j,:,:) \leftarrow
			\operatorname{\texttt{recover2d}}(\hat{I}_{\Omega_{ppx}}(q(j-1)+1,:,:),\hat{I}_D(j,:,:),j-1)$} $~~$ \Comment{Algorithm~\ref{alg:3DPPFT:recover2d}}
		
		\State{$\hat{I}_D(n+2-j,:,:) \leftarrow
			\operatorname{\texttt{recover2d}}(\hat{I}_{\Omega_{ppx}}(qn-q(j-1)+1,n+1:-1:1,n+1:-1:1),\hat{I}_D(n+2-j,:,:),j-1)$}
		
		\State{$\hat{I}_D(:,j,:) \leftarrow
			\operatorname{\texttt{recover2d}}(\hat{I}_{\Omega_{ppy}}(q(j-1)+1,:,:),\hat{I}_D(:,j,:),j-1)$} $~~$
		
		\State{$\hat{I}_D(:,n+2-j,:) \leftarrow
			\operatorname{\texttt{recover2d}}(\hat{I}_{\Omega_{ppy}}(qn-q(j-1)+1,n+1:-1:1,n+1:-1:1),\hat{I}_D(:,n+2-j,:),j-1)$}
		
		\State{$\hat{I}_D(:,:,j) \leftarrow
			\operatorname{\texttt{recover2d}}(\hat{I}_{\Omega_{ppz}}(q(j-1)+1,:,:),\hat{I}_D(:,:,j),j-1)$} $~~$
		
		\State{$\hat{I}_D(:,:,n+2-j) \leftarrow
			\operatorname{\texttt{recover2d}}(\hat{I}_{\Omega_{ppz}}(qn-q(j-1)+1,n+1:-1:1,n+1:-1:1),\hat{I}_D(:,:,n+2-j),j-1)$}
		
		\EndFor
		
		\State{$\hat{I}_D(n/2+1,n/2+1,n/2+1) \leftarrow
			\hat{I}_{\Omega_{ppx}}(3n/2+1,n/2+1,n/2+1)$} (center point)
		
		\State{${I}\leftarrow$\texttt{InvDecimatedFreq}($\hat{I}_{D}$)} \Comment{(recover $I$ from $\hat{I}_{D}$ by Algorithm~\ref{alg:3DPPFT:InvDecimated})}

	\end{algorithmic}
	\label{alg:3DPPFT:inversion}
\end{algorithm}

\begin{algorithm}
	\caption{\texttt{recover2d}. \protect This function recovers a 2D slice of the 3D pseudo-polar grid and recovers it to Cartesian grid (resample $\hat{I}$ from the set $\Omega_{c}^{(k+1)}$ to the set $\Theta^{(k)}$.)}
	\begin{algorithmic}[1]
		\REQUIRE $\hat{U}_{\Omega_{ppx}}$ - Slice of the 3D pseudo-polar Fourier transform of size $(n+1)\times(n+1)$.
		
		$\hat{U}_D$ - Already resampled slice of $\hat{I}_D$ from Algorithm~\ref{alg:3DPPFT:inversion} of size
		$(n+1)\times(n+1)$.
		
		$j$ - An integer indicating the stage within the peeling procedure.
		\ENSURE $result$ of size $(n+1) \times (n+1)$ (resampled part, $\hat{I}_D$ from Algorithm~\ref{alg:3DPPFT:inversion}) 
		
		\State $\alpha \leftarrow (n/2-j)/(n/2)$
		
		\State $m \leftarrow qn+1$
		
		\If {$j==0$}
		
		\State $result \leftarrow \hat{U}_{\Omega_{pp}}$
		
		\Return		
		
		\EndIf

		\State $y_1 \leftarrow [-n/2:n/2] \times (-2) \times q \times  \pi /m$
		
		\State $x_1 \leftarrow [-n/2:n/2] \times (-2) \times q \times \alpha \times \pi
		/m$
		
		\State $C \leftarrow \operatorname{zeros}(n+1,n+1)$
		
		\For{$k=1,\ldots,j$}
		
		\State $C(k,:) \leftarrow
		\texttt{ToeplitzResample}\left ( \hat{U}_{\Omega_{pp}}(k,1:n+1),y_1,x_1 \right )$ \Comment{Algorithm~\ref{alg:3DPPFT:toepresampling}}
		
		\State $C(n+2-k,:) \leftarrow
		\texttt{ToeplitzResample}\left (\hat{U}_{\Omega_{pp}}(n+2-k,1:n+1)y_1,x_1  \right )$
		
		\EndFor
		
		\State $x \leftarrow [-(n/2-j):(n/2-j)] \times (-2) \times q \times \pi /m$
		
		\State $y \leftarrow \left [ [-n/2:-n/2+j-1], [-n/2:n/2]\times \alpha,
		[n/2-j+1:n/2] \right ] \times (-2) \times q \times \pi/m$
		
		\State $R_1 \leftarrow \operatorname{zeros}(n+1-2j,n+1)$
		
		\For{$k=1,\ldots,n+1$}
		
		\State $R_1 \leftarrow \texttt{ToeplitzResample}\left (\left [ C(1:j,k),
		\hat{U}_{\Omega_{pp}}(:,k), C(n-j+2:n+1,k)\right ],y,x\right ) $
		
		\EndFor
		
		\State $R_2 \leftarrow \operatorname{zeros} (n+1-2j,n+1-2j)$
		
		\For{$k=1,\ldots,n+1-2j$}
		
		\State $R_2 \leftarrow \texttt{ToeplitzResample}\left ( \left [\hat{U}_D(k+j,1:j), R_1(k,:),
		\hat{U}_D(k+j,n-j+2:n+1) \right ],y,x \right ) $
		
		\EndFor
		
		\State $result \leftarrow \hat{U}_D$
		
		\State $result (1+j:n+1-j,1+j:n+1-j) \leftarrow R_2$

	\end{algorithmic}
	\label{alg:3DPPFT:recover2d}
\end{algorithm}

%
%
%
%
%
%

\subsection{Recovering $I$ from $\hat{I}$ on $\Omega_{c}$}\label{sec:invDecimatedF}
\label{invertingDecFreq}
\begin{figure}
	\centering
	\includegraphics[width=0.7\textwidth]{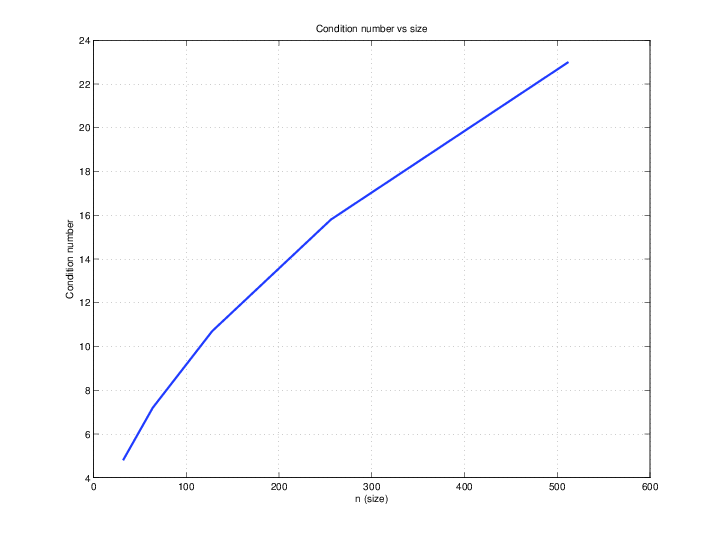}
	\caption{The maximal condition number of $F_D$ (see~\eqref{eq:FD}) obtained while applying Algorithm~\ref{alg:3DPPFT:InvDecimated} to various sizes $n$.}
	\label{fig:cond_num}
\end{figure}

The second step in Algorithm~\ref{alg:3DPPFT:inversion} (line 11) recovers the volume $I$ from the samples in $\hat{I}$ on $\Omega_{c}$ (see~\eqref{eq:DecCaGrid}). Define the operator $F_D:\mathbb{C}^n \to \mathbb{C}^{n+1}$, by
\begin{equation}\label{eq:FD}
(F_D u)_k=\sum_{j=-n/2}^{n/2-1}u(j)e^{-2\pi \imath jqk/m},\quad k=-n/2,\ldots,n/2,\quad
m=qn+1.
\end{equation}
For a volume $I$ of size $n \times n \times n$, we denote by $F_{D}^{(x)}$, $F_{D}^{(y)}$ and $F_{D}^{(z)}$ the application of $F_{D}$ to the $x$, $y$ and $z$ directions of $I$, respectively. Furthermore, we define
\begin{equation}\label{eq:Itilde}
\tilde{I} = F_{D}^{(z)} F_{D}^{(y)} F_{D}^{(x)} I.
\end{equation}
From the definitions of $F_{D}$, $\hat{I}$ and $\Omega_{c}$ in Eqs. ~\ref{eq:FD},~\ref{eq:3DPPFT:dtft3} and \ref{eq:DecCaGrid}, respectively, we conclude that $\tilde{I}$ is equal to the samples of $\hat{I}$ on the grid $\Omega_{c}$. Thus, if we apply the inverse of $F_{D}$ to the $x$, $y$, and $z$ dimensions of $\tilde{I}$ (in a separable way), we recover $I$ from the samples of $\hat{I}$ on $\Omega_{c}$. The inversion of the operator $F_D$ is described in~\cite{averbuch2008framework}. Specifically, in our case, for $m=qn+1$, the adjoint of $F_{D}$ is given by
\begin{equation*}
(F^*_Dw)_{j}=\sum_{k=-n/2}^{n/2}w_{k} e^{2\pi \imath j q k/m},\quad
j=-n/2,\ldots,n/2-1.
\end{equation*}
The matrix $F^*_{D}F_D$ is a Toeplitz matrix, whose entries are given by
\begin{equation*}
(F^*_{D}F_D)_{k,l}=\sum_{j=-n/2}^{n/2}e^{2\pi
	\imath q j(k-l)/m},\quad k,l=-n/2,\ldots,n/2-1,
\end{equation*}
and its inverse is applied as described in \eqref{eq:gohberg-semencul formula}.
The image $I$ is recovered by the application of $(F_D^*F_D)^{-1}F_D^*$ to each dimension of $\widetilde{I}$. Accurately recovering $I$ from $\tilde{I}$ requires the condition number of $F^*_{D}F_D$ to be small. The maximal condition number of $F_D$, which was obtained while applying Algorithm~\ref{alg:3DPPFT:recover2d} to various sizes $n$, is illustrated in Fig.~\ref{fig:cond_num}. The reason the term ``maximal" is used, is since during the running of the algorithm, several operators $F_D$ are used for a given $n$. As can be seen, the condition number is less than 25 even for very large volumes. The operator $F_{D}^{*}$ is applied efficiently to a vector $w$ by adding $q-1$ zeros between every two samples of $w$, applying the inverse FFT to the resulting vector, and keeping the $n$ central elements of the FFT-transformed vector (see Algorithm~\ref{alg:3DPPFT:adjF3D} for a detailed description). The pseudo-code of the algorithm for recovering the volume from its samples in $\Omega_{c}$ is described in Algorithm~\ref{alg:3DPPFT:InvDecimated}.

\begin{algorithm}
	\caption{\texttt{AdjFDecimated}: Applying $F_D^*$ of~\eqref{eq:FD} to a vector}
	\begin{algorithmic}[1]
		\REQUIRE $y$ - Input vector of odd length.\\
		$q$ - integer (decimation factor).
		\ENSURE $x=F_D^*y$
		
		\State $n \leftarrow$length$(y)-1$
		\State $m \leftarrow qn+1$
		\State $z \leftarrow\operatorname{zeros}(m,1)$
		\State $z(1:q:m) \leftarrow y$
		\State $l \leftarrow$length$(z)$
		\State $p_f \leftarrow$ floor$(l/2)$
		\State $p_c \leftarrow$ ceil$(l/2)$
		\State $x \leftarrow z([p_f+1:l,~1:p_f])$
		\State $x_{IFFT} \leftarrow$IFFT$(x)$
		\State $x \leftarrow m\cdot x_{IFFT}([p_c+1:l,~1:p_c])$
		\State $x \leftarrow x(n+1:2n)$
	\end{algorithmic}
	\label{alg:3DPPFT:adjF3D}
\end{algorithm}

\begin{algorithm}
	\caption{\texttt{InvDecimatedFreq}. \protect This procedure inverts the decimated frequencies and restores $I$ from $\tilde{I}$ of Eq. \eqref{eq:Itilde}}
	\begin{algorithmic}[1]
		\REQUIRE $\tilde{I}$ - Matrix of size $(n+1) \times (n+1) \times (n+1)$.
		\ENSURE $I$ - Volume of size $n \times n \times n$ such that~\eqref{eq:Itilde} holds.
		
		\State $c \leftarrow$ zeros$(n,1)$
		\State $m \leftarrow qn+1$
		\For{ $k=-n/2,\ldots, n/2-1$}
		\For{ $l=-n/2,\ldots, n/2$}
		\State	$c(k+n/2+1) \leftarrow c(k+n/2+1)+e^{q\pi \imath l(-n/2-k)/m}$
		\EndFor
		\EndFor
		\State $(M_1,M_2,M_3,M_4) \leftarrow$ \texttt{ToeplitzInv}($c$) \Comment{Algorithm~\ref{alg:3DPPFT:topinv}}
		\State $D_i \leftarrow$ \texttt{ToeplitzDiag}($M_i$), ~$i=1,\ldots,4$ \Comment{Algorithm~\ref{alg:3DPPFT:topdiag}}
		\State $I_{1} \leftarrow\operatorname{zeros}(n+1,n+1,n)$
		\State $I_{2} \leftarrow\operatorname{zeros}(n+1,n,n)$
		\State $I \leftarrow\operatorname{zeros}(n,n,n)$
		\For{ $k=1,\ldots, n+1$}
		\For{ $l=1, \ldots, n+1$}
		\State $v \leftarrow \widetilde{I}(k,l,:)$
		\State $v \leftarrow $\texttt{AdjFDecimated}$(v)$ \Comment{Algorithm~\ref{alg:3DPPFT:adjF3D}}
		\State $u \leftarrow$\texttt{ToeplitzInvMul}$(D_1,D_2,D_3,D_4,v)$ \Comment{Algorithm~\ref{alg:3DPPFT:topinvmul}}
		\State $I_{1}(k,l,:) \leftarrow u$
		\EndFor
		\EndFor
		\For{ $k=1,\ldots, n+1$}
		\For{ $l=1,\ldots, n$}
		\State $v \leftarrow I_{1}(k,:,l)$
		\State $v \leftarrow $\texttt{AdjFDecimated}$(v)$
		\State $u \leftarrow$\texttt{ToeplitzInvMul}$(D_1,D_2,D_3,D_4,v)$
		\State $I_{2}(k,l,:) \leftarrow u$
		\EndFor
		\EndFor	
		\For{ $k=1,\ldots, n$}
		\For{ $l=1,\ldots, n$}
		\State $v \leftarrow I_{2}(:,k,l)$
		\State $v \leftarrow $\texttt{AdjFDecimated}$(v)$
		\State $u \leftarrow$\texttt{ToeplitzInvMul}$(D_1,D_2,D_3,D_4,v)$
		\State $I(:,k,l) \leftarrow u$
		\EndFor
		\EndFor			
	\end{algorithmic}
	\label{alg:3DPPFT:InvDecimated}
\end{algorithm}

\subsection{Complexity Analysis}\label{sec:complexity}
The inversion Algorithm~\ref{alg:3DPPFT:inversion} consists of two steps: Resampling (lines 1-12) and recovering $I$ from $\hat{I}$ on $\Omega_c$ (line 13). The first step, resampling, which utilizes Toeplitz matrices (Algorithm~\ref{alg:3DPPFT:toepresampling}), takes $\mathcal{O}(n\text{log}n)$ operations when preprocessing took place. This function (Toepliz-based resampling) is called $\mathcal{O}(n)$ times in the $\operatorname{\texttt{recover2d}}$ function (Algorithm~\ref{alg:3DPPFT:recover2d}), which in turn is called also $\mathcal{O}(n)$ times by Algorithm~\ref{alg:3DPPFT:inversion}. A total of $\mathcal{O}(n^3\text{log}n)$ operations are needed for the resampling step.

The second step, recovering $I$ from $\hat{I}$ on $\Omega_c$, takes  $\mathcal{O}(n^3\text{log}n)$ operations. This yields a total computational complexity of $\mathcal{O}(n^3\text{log}n)$ operations. If the preprocessing is done in run time, the total complexity becomes $\mathcal{O}(n^4)$ operations, due to the Durbin-Levinson algorithm used for inverting Toeplitz matrices. The computational complexity can be reduced from $\mathcal{O}(n^4)$ operations to $\mathcal{O}(n^3\text{log}^2 n)$ by solving~\eqref{eq:xyequation} as  described in~\cite{kailath1989generalized}.

The total storage requirement for $n$ applications of Algorithm~\ref{alg:3DPPFT:recover2d} is $4n^2$. This storage is used to store the diagonal Toeplitz matrices $D_1,\ldots,D_4$.

\section{Numerical results}
\label{numerical}
We implemented Algorithm~\ref{alg:3DPPFT:inversion} in MATLAB and applied it to several volumes of sizes $n \times n \times n$ where $n=32,64,...,512$. All the experiments were executed on a Linux machine with two Intel Xeon processors (CPU X5560) running at 2.8GHz, with 8 cores in total and $96$GB of RAM. All the 3D experiments were performed with $q=3$ (see Eq.~\ref{eq2.4}).

\begin{figure}
	\centering
	\includegraphics[width=0.7\textwidth]{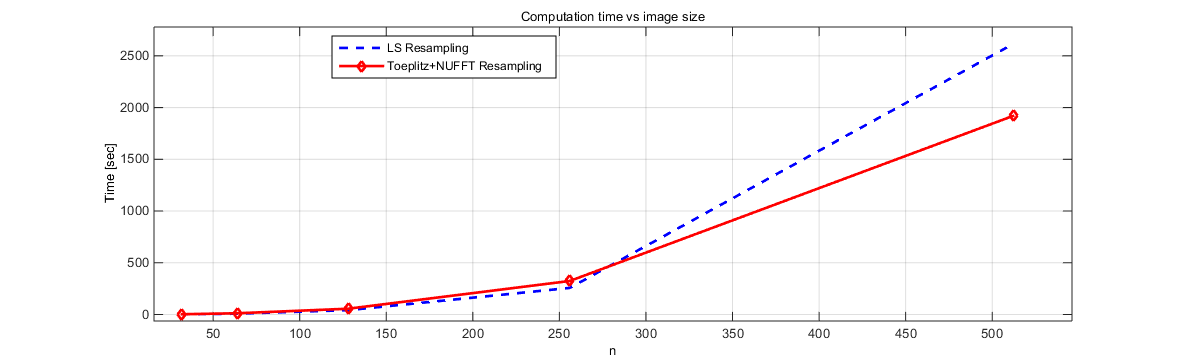}
	\caption{Runtime comparison between LS-based resampling and Toeplitz-NUFFT-based resampling (Algorithm~\ref{alg:3DPPFT:toepresampling}) for volumes of sizes $n \times n \times n$, consisting of PPFT applied to a normally distributed i.i.d. random samples with zero mean and unit variance.}
	\label{fig:3dppft_results random}
\end{figure}

\begin{figure}
	\centering
	\includegraphics[width=0.7\textwidth]{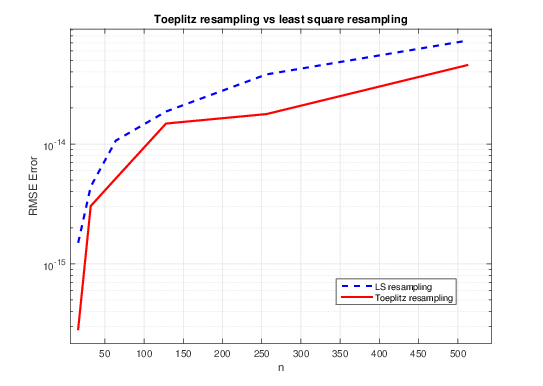}
	\caption{Error comparison between LS-based resampling and Toeplitz-NUFFT-based resampling (Algorithm~\ref{alg:3DPPFT:toepresampling}) for  volumes of sizes $n\times n \times n$, consisting of PPFT applied to a normally distributed random samples with zero mean and unit variance.}
	\label{fig:compare3d_err}
\end{figure}

Algorithm~\ref{alg:3DPPFT:recover2d} is based on a series of one-dimensional resampling operations. We compare two methods for implementing this one-dimensional resampling -- least-squares-based (LS) approach and Toeplitz-NUFFT-based approach, both described in Section~\ref{sec:1d resampling algorithm}. The LS-based approach, described in Eq.~\ref{eq:resamplingLS}, consists of finding the coefficients vector $\alpha$ of a trigonometric polynomial, followed by direct evaluation of the polynomial at the resampling points. Using this resampling approach in Algorithm~\ref{alg:3DPPFT:recover2d} results in computational complexity of $\mathcal{O}(n^3)$ operations (excluding a preprocessing step). Nevertheless, its implementation in MATLAB is highly optimized. The Toeplitz-NUFFT-based approach, described in Section~\ref{sec:1d resampling algorithm} and Algorithm~\ref{alg:3DPPFT:toepresampling}, results in computational complexity of Algorithm~\ref{alg:3DPPFT:recover2d} of $\mathcal{O}(n^2 \text{log}n + n^{2} \log 1/\epsilon)$ operations. The two resampling methods are compared by using them as the underlying one-dimensional resampling in Algorithm~\ref{alg:3DPPFT:recover2d}. Both algorithms use a preprocessing step which is excluded from the reported running time. The error incurred by the two algorithms is measured by
\begin{equation*}
E(I,I_r) = \sqrt{\frac{\sum_{k,l,m}\abs{I[k,l,m]-I_r[k,l,m]}^2}{\sum_{k,l,m}\abs{I[k,l,m]}^2}},
\end{equation*}
where $I[k,m,n]$ is the original volume and $I_r[k,m,n]$ is the reconstructed volume.
Running times for both methods for volumes of sizes $n \times n \times n$, where $n=32,64,...,512$, are shown in Fig.~\ref{fig:3dppft_results random}. The volumes in this experiment consist of PPFT transformed random independent samples from a  normal distribution with zero mean and unit variance. Figure \ref{fig:compare3d_err} compares between the reconstruction errors of Algorithm \ref{alg:3DPPFT:toepresampling} and those of LS-based resampling. The results show that both methods are comparable in terms of accuracy.

\begin{figure}
	\centering
	\includegraphics[width=0.7\textwidth]{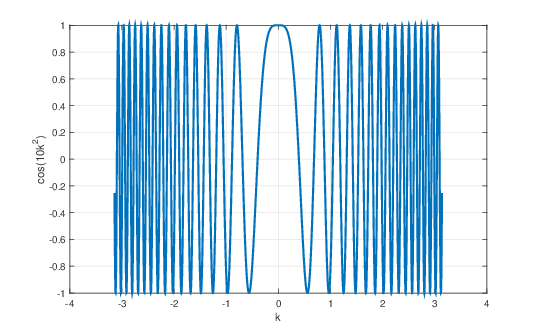}
	\caption{One-dimensional Chirp signal, used for testing the LS-based resampling and the Toeplitz-NUFFT-based resampling.}
	\label{fig:chirp1d}
\end{figure}

\begin{figure}
	\centering
	\includegraphics[width=0.7\textwidth]{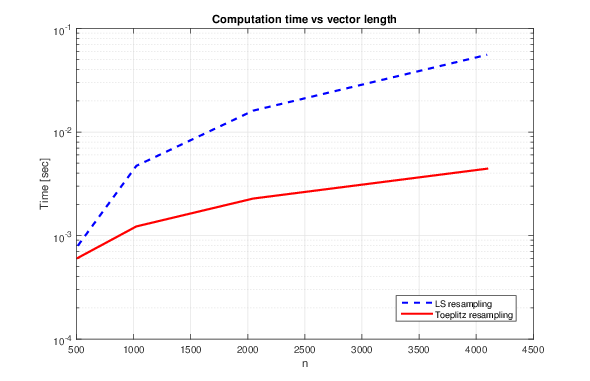}
	\caption{Runtime comparison between the LS-based resampling and the Toeplitz-NUFFT-based resampling (Algorithm~\ref{alg:3DPPFT:toepresampling}) for a one-dimensional Chirp signal of length $n$.}
	\label{fig:compare1d}
\end{figure}

\begin{figure}
	\centering
	\includegraphics[width=0.7\textwidth]{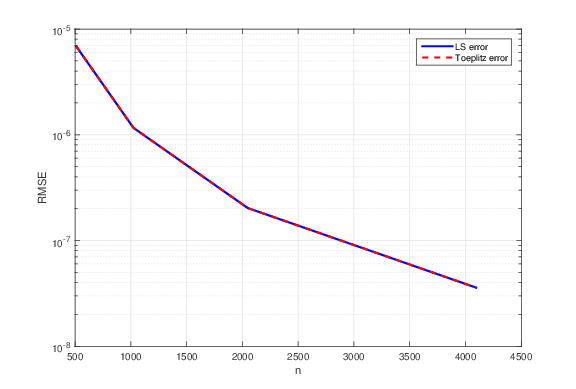}
	\caption{Error comparison between the LS-based resampling and Toeplitz-NUFFT-based resampling (Algorithm~\ref{alg:3DPPFT:toepresampling}) for a one-dimensional chirp signal of length $n$.}
	\label{fig:compare1derr}
\end{figure}

Since it is currently impossible to process 3D volumes of sizes larger than $n=512$, we compare between the LS-based resampling and the Toeplitz-based resampling by applying both one-dimensional resampling methods to a one-dimensional Chirp signal defined by $\cos(10k_{j}^2)$, $k_{j}=-\pi+\frac{2\pi}{n}j$,  $j=0,\ldots,n$, $n=512,\ 1024,\ 2048,\ 4096$.  The signal's samples given at $k_j$ were resampled to $0.3k_j$ using LS-based resampling and Toeplitz-based resampling. The one-dimensional Chirp signal is displayed in Fig.~\ref{fig:chirp1d}. The running time of the two methods, which does not include pre-processing timing, appears in Fig.~\ref{fig:compare1d}. The approximation errors of both methods are practically identical as shown in Fig.~\ref{fig:compare1derr}.

\begin{figure}
	\centering
	\includegraphics[width=0.7\textwidth]{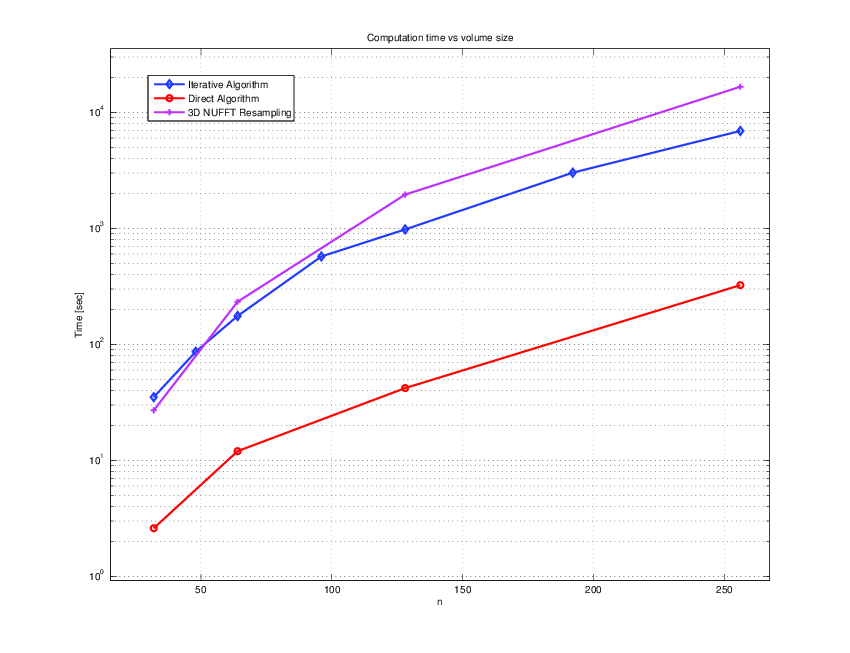}
	\caption{Comparison between the running time of Algorithm~\ref{alg:3DPPFT:inversion}, the iterative algorithm~\cite{iterativeradon2014}, and the NUFFT-based algorithm~\cite{fenn2007computation}.}
	\label{fig:comparison_old_new}
\end{figure}

Next, we compare between the direct inverse PPFT in Algorithm~\ref{alg:3DPPFT:inversion}, the iterative algorithm described in~\cite{iterativeradon2014}, and a single iteration of an implementation that computes the Gram operator of the PPFT using the NUFFT as suggested by~\cite{fenn2007computation}. For the latter method, several iterations are required, whose number depends on the condition number of the operator. However, to keep all timings within the same scale, we compare the other two algorithms to only a single iteration of the latter one. The results are illustrated in Fig.~\ref{fig:comparison_old_new}, which shows that Algorithm~\ref{alg:3DPPFT:inversion} is faster than the iterative algorithm~\cite{iterativeradon2014} as well as than a single iteration of the 3D NUFFT-based algorithm proposed by~\cite{fenn2007computation}.  Results for $n=512$ do not appear in Fig.~\ref{fig:comparison_old_new} (unlike Fig.~\ref{fig:3dppft_results random}) as it was impossible to process volumes of that size using~\cite{iterativeradon2014,fenn2007computation}.

Given the long running time of the 3D NUFFT-based algorithm, we executed above only a single iteration of this algorithm. However, to compare accuracy of our proposed method with that of the 3D NUFFT-based algorithm, we execute next the 3D NUFFT-based algorithm until the error becomes smaller than $10^{-12}$ or the number of iterations exceeds 100. Due to time constraints, we used only $n=16,\ 32,\ 64$. The input for each case was of size $n\times n \times n$ of random normally distributed i.i.d. samples with zero mean and unit variance. We implemented the 3D NUFFT-based algorithm with the same preconditioner as in \cite{iterativeradon2014}. For $n=16$ the algorithm took 112 seconds and the resulting error was $7.25 \times 10^{-13}$, obtained after 40 iterations; for $n=32$, it took $2,710$ seconds with error of $7.5\times 10^{-13}$ after 100 iterations; for $n=64$ it took $22,477$ seconds (more than 6 hours) with error of $1.61\times 10^{-12}$ after 100 iterations. The iterative algorithm described in~\cite{iterativeradon2014} gives similar errors to the direct inversion algorithm (Algorithm \ref{alg:3DPPFT:inversion}) whose error appears in Fig.~\ref{fig:compare3d_err}.

\begin{figure}
	\centering
	\subfloat[Dolphin, $n=64$]{\label{dolphin}
		\includegraphics[width=0.3\textwidth]{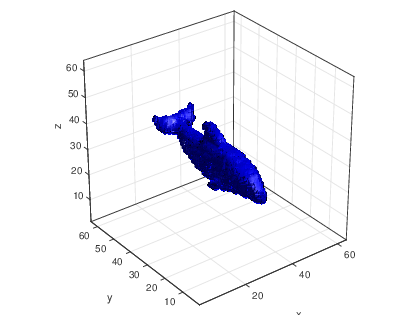}}
	\subfloat[Bird, $n=128$]{\label{bird}
		\includegraphics[width=0.3\textwidth]{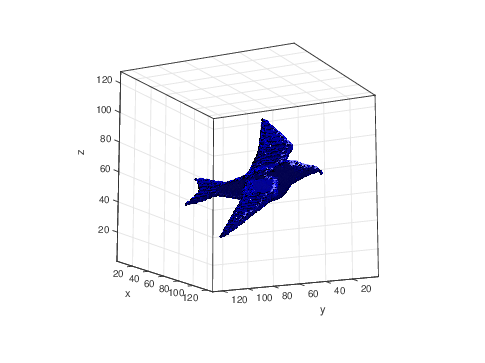}}
	\subfloat[Cup, $n=256$]{\label{cup}
		\includegraphics[width=0.3\textwidth]{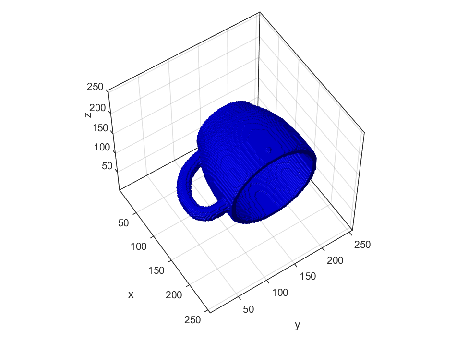}}
	\caption{Test volumes.}
	\label{fig:realvols}
\end{figure}

Finally, we tested the algorithm on real volumes of different sizes: a dolphin of size $64 \times 64 \times 64$ (Fig.~\ref{dolphin}), a bird of size $128 \times 128 \times 128$ (Fig.~\ref{bird}), and a 3D cup of size $256\times 256 \times 256$ (Fig.~\ref{cup}). The volumes were taken from the McGill 3D shape dataset \cite{siddiqi2008retrieving} and PPFT was applied to each one them to be used as an input for the inversion algorithm.
The results appear in Table \ref{tab:realshapes}.

\begin{table}
	\centering
	\begin{tabular}{c*{6}{c}c}
		Name & Volume & Time [sec]& RMSE \\
		\hline
		Dolphin & $64\times 64 \times 64$ & 7.43 & $1.69 \times 10^{-15}$   \\
		Bird            & $128\times 128 \times 128$ & 40.15 & $3.6\times 10^{-15}$  \\
		Cup           & $256\times 256 \times 256$ & 270 & $1.25 \times 10^{-14}$   \\
	\end{tabular}
	\caption{Results for volumes of Fig.~\protect\ref{fig:realvols}.}
	\label{tab:realshapes}
\end{table}

\section{Conclusions}\label{sec:conclusions}
In this paper, a new algorithm for inverting the 3D pseudo-polar Fourier transform (PPFT) is described. The algorithm processes at each iteration a two-dimensional slice of the input, where each such processing uses only one-dimensional operations. The main component of the algorithm is fast resampling of univariate trigonometric polynomials. The resampling is implemented using a 1D non-uniform Fourier transform together with fast algorithms for Toeplitz matrices.
The algorithm is not iterative, and requires a fixed amount of time that depends only on the size of the input. Moreover, the algorithm has low memory requirements, allowing to process large 3D datasets in a reasonable time. The performance of the algorithm is demonstrated on volumes as large as $512 \times 512 \times 512$ in double precision.

\begin{appendices}
\section{Toeplitz solvers algorithms}
\begin{algorithm}
	\caption{\texttt{ToeplitzResample}: Fast resampling of a univariate trigonometric polynomial.}
	\begin{algorithmic}[1]
		\REQUIRE $y$ - vector of length $n$.\\
        $f$ - values of the trigonometric polynomial $\hat{I}$ on the set $y$.\\
		$x$ - vector of length $n$.
		\ENSURE $g$ - values of $\hat{I}$ on the set $x$.
		
		\State $k \leftarrow [-n/2:n/2-1]$ \Comment{$k$ is a row vector}
		
		\State $c \leftarrow \text{length}(y)\times \text{NUFFT}_1(y, e^{\imath  yk(1)},-1, n)$  \Comment{$k(1)$ is the first element of $k$}
		\State $(M_1,M_2,M_3,M_4 ) \leftarrow$ \texttt{ToeplitzInv}$(c,c)$ \Comment{Can be computed in the preprocessing step}
        \For {$i=1,\ldots,4$}
    		\State $D_i \leftarrow$\texttt{ToeplitzDiag}($M_i(:,1), M_i(1,:)$) \Comment{Can be computed in the preprocessing step}
        \EndFor
		\State $v\leftarrow $ \text{length} $(y)\times$ NUFFT$_1(y,f,-1,n)$
		\State $\alpha \leftarrow $\texttt{ToeplitzInvMul}$(D_1,D_2,D_3,D_4,v)$
		\State $g \leftarrow $ NUFFT$_2$ $(x,1,n,\alpha)$
		
	\end{algorithmic}
	\label{alg:3DPPFT:toepresampling}
\end{algorithm}
NUFFT$_k$ ($k=1,2$) refers to the type of the NUFFT - see Eqs. ~\ref{eq:nuff1} and~\ref{eq:nuff2}. The parameters $-1$ and $1$ in lines $6,11,13$ in Algorithm~\ref{alg:3DPPFT:toepresampling} refer to the sign of $\imath$ in the complex exponent. The function \texttt{ToeplitzInv} is described in Algorithm~\ref{alg:3DPPFT:topinv} and computes the Gohberg-Semencul factorization of a symmetric Toeplitz matrix.

\begin{algorithm}
	\caption{\texttt{ToeplitzInv}: Factorize a symmetric Toeplitz matrix using the Gohbcerg-Semencul formula.}
	\begin{algorithmic}[1]
		\REQUIRE $c$ - First column (row) of a symmetric Toeplitz matrix.
		\ENSURE $M_1,M_2,M_3,M_4$ - Gohberg-Semencul factorization.
		\State $n \leftarrow$length$(c)$
		\State $e_0 \leftarrow\operatorname{zeros}(n,1)$
		\State $e_n \leftarrow\operatorname{zeros}(n,1)$
		\State $e_0(1) \leftarrow 1$
		\State $e_n(n) \leftarrow 1$
		\State Solve $T_n(c,c)x=e_0$
		\State Solve $T_n(c,c)y=e_n$
		\State Compute $M_1$ using~\eqref{ppft:m1eq}
		\State Compute $M_2$ using~\eqref{ppft:m2eq}
		\State Compute $M_3$ using~\eqref{ppft:m3eq}
		\State Compute $M_4$ using~\eqref{ppft:m4eq}
	\end{algorithmic}
	\label{alg:3DPPFT:topinv}
\end{algorithm}

\begin{algorithm}
	\caption{\texttt{ToeplitzDiag}: Diagonalize the circular embedding of a Toeplitz matrix}
	\begin{algorithmic}[1]
		\REQUIRE $c, r$ - First column/row of the Toeplitz matrix.
		\ENSURE $D$ - Diagonal form of $T_{n}(c,r)$ embedded in a circulant matrix.
		
		\State $n \leftarrow$ length$(c)$
		\State $C \leftarrow [c; 0; r(n:-1:2)]$
		\State $D \leftarrow$ FFT$(C)$
	\end{algorithmic}
	\label{alg:3DPPFT:topdiag}
\end{algorithm}

\begin{algorithm}
	\caption{\texttt{ToeplitzInvMul}: Mulitply an inverse Toeplitz matrix $A^{-1}$ of size $n \times n$ by a vector $v$, given the diagonal forms of the
		Gohberg-Semencul factors of $A^{-1}$  in $\mathcal{O}(n \text{log} n)$ operations.}
	\label{alg:3DPPFT:topinvmul}
	\begin{algorithmic}
		\REQUIRE $D_1,D_2,D_3,D_4$ - Diagonal forms of the Gohberg-Semencul factors of the inverse Toeplitz matrix. \\
		$v$ - vector of length $n$.
		\ENSURE $u$ - The product $A^{-1}v$.
		\State $u_1 \leftarrow$\texttt{ToeplitzMul}$(D_2,v)$
		\State $u_1 \leftarrow$\texttt{ToeplitzMul}$(D_1,u_1)$
		\State $u_2 \leftarrow$\texttt{ToeplitzMul}$(D_4,v)$
		\State $u_2 \leftarrow$\texttt{ToeplitzMul}$(D_3,u_2)$
		\State $u \leftarrow u_1-u_2$
	\end{algorithmic}
\end{algorithm}

\begin{algorithm}
	\caption{\texttt{ToeplitzMul}: Multiply a Toeplitz matrix $A$ of size $n \times n$, whose diagonal factor is $D$, by a vector $v$ in $\mathcal{O}(n \text{log}n)$ operations.}
	\begin{algorithmic}[1]
		\REQUIRE $D$ - Diagonal factor of $A$ (computed using Algorithm~\ref{alg:3DPPFT:topdiag}). \\
         $v$ - vector of length $n$.
		\ENSURE $u$ - the product $A v$.
		
		\State $n \leftarrow$ length$(v)$
		\State $v \leftarrow [v;\text{zeros}(n,1)]$
		\State $w \leftarrow \text{FFT}(v)$
		\State $u \leftarrow \text{IFFT}(D.*w)$
		\State $u \leftarrow u(1:n)$
	\end{algorithmic}
	\label{alg:3DPPFT:topmul}
\end{algorithm}

\end{appendices}
\newpage
\bibliographystyle{plain}
\bibliography{3d_ppft_new_ideas}

\end{document}